\documentclass[11pt]{amsart}

\usepackage{amsfonts,amsmath,amssymb}

{\theoremstyle{plain}
\newtheorem{theorem}{Theorem}    
}

\newcommand{\Cc}{\mathbb{C}}  
\newcommand{\Pp}{\mathbb{P}}

\newcommand{\defi}[1]{\emph{#1}}

\renewcommand{\epsilon}{\varepsilon}

\renewcommand{\ge}{\geqslant}

\newcommand{\grad}{\mathop{\mathrm{grad}}\nolimits}
\newcommand{\Baff}{{\mathcal{B}_\mathit{\!aff}}}
\newcommand{\Binf}{{\mathcal{B}_\infty}}
\newcommand{\B}{{\mathcal{B}}}

\newcommand{\pr}{\mathrm{pr}}
\newcommand{\aff}{\mathit{aff}}

\newcommand{\PP}{{\mathcal{P}}}
\newcommand{\CC}{{\mathcal{C}}}
\newcommand{\mm}{{\mathfrak{m}}}

\newcommand{\ttt}[1]{\texttt{#1}}
\newcommand{\Singular}{\textsc{Singular}}

\begin{document}
\title{Computation of Milnor numbers and critical values at infinity}
\author{Arnaud Bodin}
\date{\today}
\address{laboratoire AGAT, UFR de Math\'ematiques, Universit\'e Lille I, 59655 Villeneuve d'Ascq}
\email{Arnaud.Bodin@agat.univ-lille1.fr \\
http://www-gat.univ-lille1.fr/\~{}bodin
}

\begin{abstract}
We describe how to compute
topological objects associated to a complex polynomial map of $n
\ge 2$ variables with isolated singularities. These objects are:
the affine critical values, the affine Milnor numbers for all irregular fibers, the
critical values at infinity, and the Milnor numbers at infinity for all irregular fibers.
Then for a family of polynomials we detect parameters where the topology of the polynomials 
can change. Implementation and examples are given with the 
computer algebra system \Singular.
\end{abstract}

\maketitle

\section{Introduction}

\subsection{Review on the local case}

Let $g : \Cc^n,0 \longrightarrow \Cc,0$ be a germ of polynomial map with 
isolated singularities. One of the most important topological object attached to $g$
is its \defi{local Milnor number} \cite{Mi}:
$$\mu_0 = \dim_\Cc {\Cc\{x_1,\ldots,x_n\}} / {\mathrm{Jac}(g)}$$
where $\mathrm{Jac}(g)= (\frac{\partial g}{\partial x_1},\ldots,\frac{\partial g}{\partial x_n})$
is the Jacobian ideal of $g$.
It is possible to compute $\mu_0$ with the help of a Gr\"{o}bner base. For example such a computation
motivates the computer algebra system \Singular, \cite{sing}.

Now we consider a family $(g_s)_{s\in[0,1]}$, with  $g_s : \Cc^n,0 \longrightarrow \Cc,0$
germs of isolated singularities, such that $g_s$ is  a smooth function of $s$. To each 
$s\in [0,1]$ we associate the local Milnor number $\mu_0(g_s)$. The main topological result 
for families is L\^e-Ramanujam-Timourian $\mu$-constant theorem.

\begin{theorem}[\cite{LR,Tim}]
\label{th:mucstloc}
If $n\not=3$ and $\mu_0(g_s)$ is constant ($s\in [0,1]$) then the family
$(g_s)_{s\in[0,1]}$ is a topologically trivial family.
\end{theorem}

\subsection{Motivation and aims for the global case}

Now we consider a polynomial function $f : \Cc^n \longrightarrow \Cc$. The study of the topology
of $f$ is not just the glueing of local studies because of the behaviour of $f$ at infinity, 
see \cite{Br}.
To the polynomial $f$ we attach ``Milnor numbers'' $\mu$, $\lambda$ and finite sets of critical values
 $\Baff$, $\Binf$, $\B= \Baff \cup \Binf$
(see the definitions below). 
The first aim of this work is to compute these objects and to give the topology
of the fibers $f^{-1}(c)$ for all $c\in \Cc$.

There is a global version of the local $\mu$-constant theorem (see Theorem \ref{th:mucst})
where the Milnor number $\mu_0$ is replaced by a \defi{Milnor multi-integer} 
$\mm =  (\mu,\#\Baff,\lambda, \#\Binf, \#\B)$.
In order to verify if $\mm(f_s)$ remains constant in a family $(f_s)_{s\in[0,1]}$
it is not possible to compute $\mm(f_s)$ for infinitely many values.
The second aim of the work is to give (and compute) a finite set $\mathcal{S}'$
such that $\mm(f_s)$ is constant for $s \in [0,1]\setminus \mathcal{S}'$.

\medskip

The rest of this section is devoted to the definitions and the results.

\subsection{Critical values}

Let $f : \Cc^n \longrightarrow \Cc$ be a polynomial map, $n\ge2$.
By a result of Thom \cite{T} there is a minimal
 \defi{set of critical values} $\B$ of point of $\Cc$ such that
$f : f^{-1}(\Cc\setminus \B) \longrightarrow \Cc\setminus \B$ is a
fibration.

\subsection{Affine singularities}
We suppose that \defi{affine singularities are isolated}
\emph{i.e.} that the set $\{ x \in \Cc^n \ | \ \grad_f x = 0\}$ is
a finite set. Let $\mu_c$ be the sum of the local Milnor numbers
at the points of $f^{-1}( c )$. Let
$$\Baff = \big\{ c \ | \ \mu_c > 0\big\} \quad \text { and }\quad \mu = \sum_{c\in\Cc} \mu_c$$
be the \defi{affine critical values} and the \defi{affine Milnor
number}.

\subsection{Singularities at infinity}

See \cite{Br}. Let $d$ be the degree of $f : \Cc^n \longrightarrow
\Cc$, let $f = f^d+f^{d-1}+\cdots+f^0$ where $f^j$ is homogeneous
of degree $j$. Let $\bar f (x,z)$ (with $x = (x_1,\ldots,x_n)$)
be the homogenisation of $f$ with the new variable $z$: $\bar
f(x,z) = f^d(x)+f^{d-1}(x)z+\ldots +f^0(x)z^d$. Let
$$X = \left\lbrace ((x:z),t)\in \Pp^n\times \Cc \ | \ \bar
f(x,z)-cz^d=0 \right\rbrace.$$
Let $\mathcal{H}_\infty$ be the hyperplane at infinity of $\Pp^n$
defined by $(z=0)$. The singular locus of $X$ has the form
$\Sigma \times \Cc$ where 
$$\Sigma = \left\lbrace (x:0) \ | \
\frac{\partial f^d}{\partial x_1}=\cdots =\frac{\partial
f^d}{\partial x_n}= f^{d-1} =0 \right\rbrace \subset
\mathcal{H}_\infty.$$ 
We suppose that $f$ has
\defi{isolated singularities at infinity} that is to say that $\Sigma$ is
finite. This is always true for $n=2$. 
We say that $f$ has \defi{strong isolated singularities at infinity} if 
$$\Sigma' = \left\lbrace (x:0) \ | \
\frac{\partial f^d}{\partial x_1}=\cdots =\frac{\partial
f^d}{\partial x_n}=0 \right\rbrace$$ 
is finite.

For a point $(x:0) \in
\mathcal{H}_\infty$, assume, for example, that $x =
(x_1,\ldots,x_{n-1},1)$ and set $\check{x} = (x_1,\ldots,x_{n-1})$
and
$$F_c(\check{x},z) = \bar f(x_1,\ldots,x_{n-1},1)-cz^d.$$ Let
$\mu_{\check{x}}(F_c)$ be the local Milnor number of $F_c$ at the
point $(\check{x},0)$. If $(x:0)\in \Sigma$ then
$\mu_{\check{x}}(F_c)>0$. For a generic $s$,
$\mu_{\check{x}}(F_s) =  \nu_{\check{x}}$, and for finitely many
$c$, $\mu_{\check{x}}(F_c) > \nu_{\check{x}}$. We set
$\lambda_{c,\check{x}} = \mu_{\check{x}}(F_c) - \nu_{\check{x}}$,
$\lambda_c = \sum_{(x:0)\in\Sigma} \lambda_{c,\check{x}}$. Let
$$\Binf = \big\{ c\in\Cc \ | \  \lambda_c > 0\big\} \quad \text {
and } \quad \lambda = \sum_{c\in\Cc} \lambda_c$$ be the
\defi{critical values at infinity} and the \defi{Milnor number at
infinity}.

We can now describe the set of critical values $\B$ as
follows (see \cite{HL} and \cite{Pa}):
$$\B = \Baff \cup \Binf.$$ Moreover by \cite{HL} and \cite{ST} for
all $c \in \Cc$, $f^{-1}(c)$ has the homotopy type of a wedge of
$\mu+\lambda-\mu_c-\lambda_c$ spheres of real dimension $n-1$.

\subsection{Families of polynomials}

To a polynomial we associate its \defi{Milnor multi-integer} 
$\mm = (\mu,\#\Baff,\lambda, \#\Binf, \#\B)$.
Two polynomials maps $f,g : \Cc^n \longrightarrow \Cc$ are \defi{topologically equivalent}
if there exist homeomorphisms $\Phi : \Cc^n \longrightarrow \Cc^n$ and $\Psi :  \Cc
 \longrightarrow \Cc$ 
such that $f\circ \Phi = \Psi \circ g$.
The Milnor multi-integer is a topological invariant, that is to say if $f$ and $g$ are topologically equivalent
then $\mm(f) = \mm(g)$.
We recall a result of \cite{Bo, BT} that is kind of converse of this property.

Let $(f_s)_{s\in[0,1]}$ be a family of polynomials, such that $f_s$ has \emph{strong}
isolated singularities at infinity and isolated affine singularities for all $s\in[0,1]$.
For each $s\in [0,1]$ we consider the Milnor multi-integer of $f_s$,
$\mm(f_s) = (\mu(s),\#\Baff(s),\lambda(s), \#\Binf(s), \#\B(s))$.
We suppose that the coefficients of the family are polynomials in $s$ and
that the degree $\deg f_s$ is constant.
\begin{theorem}[\cite{Bo,BT}]
\label{th:mucst}Let $n\not=3$.
If $\mm(f_s)$ is constant ($s\in [0,1]$),
then $f_0$ is topologically equivalent to $f_1$.
\end{theorem}

How to verify the hypotheses from a computable point of view ?
It is not possible to compute $\mm(f_s)$ for infinitely many $s\in [0,1]$.
But in fact $\mm(f_s)$ is constant except for finitely many $s$, we denote by
$\mathcal{S}$ the set of these \defi{critical parameters}.

In paragraph \ref{sec:fam} we give a computation of a finite set $\mathcal{S}'$
such that 
$$\mathcal{S} \subset \mathcal{S}'.$$
Now to check if a value $s \in \mathcal{S}'$ is in $\mathcal{S}$
we compute $\mm(f_s)$ and we compare it with $\mm(f_{s'})$ where $s'$ is any value of 
$[0,1] \setminus \mathcal{S}'$;
 now  $s \in \mathcal{S}$ if and only if $\mm(f_s) \not= \mm(f_{s'})$.

\subsection{Implementation}

The results of this paper have been implemented in 
 two libraries \ttt{critic} and \ttt{defpol}. The first one enables to calculate
all the objects defined above: $\Baff$, $\mu$, $\mu_c$ for
$c\in \Baff$ ; $\Binf$, $\lambda$, $\lambda_c$ for $c\in\Binf$.
These programs are written for \Singular, \cite{sing}. It is based
on polar curves and on the article of D.~Siersma and M.~Tib\u ar,
\cite{ST}.
For polynomials in two variables ($n=2$) a program in \textsc{Maple}
has been written by G.~Bailly-Ma\^{\i}tre, \cite{BM}, based on
a discriminant formula of H\`a~H.V., \cite{H}.
For families of polynomials the second library computes a finite set $\mathcal{S}'$ 
that contains the critical parameters.

This research has partially been supported by a Marie Curie Individual Fellowship
of the European Community (HPMF-CT-2001-01246).

\section{Milnor numbers and critical values in affine space}
\label{sec:aff}

\subsection{Milnor number}

The computation of the affine Milnor number $\mu$ is easy and
well-known (see \cite{sing} for example). Let $f \in \Cc[x_1,\ldots,x_n]$.
Let $J$ be the Jacobian
ideal of the partial derivative $(\partial f/\partial x_i)_i$. Then
$\mu$ is the vector space dimension (over $\Cc$) of a Gr\"{o}bner basis
 of the quotient $\Cc[x_1,\ldots,x_n] / J$.

\subsection{Critical values}

We add a new variable $t$. We consider the variety
$$ C = \left\lbrace (x,t)\in \Cc^n\times \Cc \ | \
f(x)-t=0 \ \text{ and }\  \grad_f x = 0\right\rbrace.$$
The critical values are the projection of $C$ on the $t$-coordinate:
$\Baff = \pr_t(C)$.

\subsection{Milnor number of a fiber}

Set $c\in \Cc$. We would like to compute $\mu_c$ the sum of the Milnor numbers
of the points of $f^{-1}(c)$. Let $J$ be the Jacobian ideal of $f$ and set
$x$ a critical point. We denote by $J_x$ the localization of $J$ at $x$.
Let $I_x = (t-c,J_x)$, the dimension of $I_x$ is equal to the Milnor
number of $f$ at $x$. For $k \ge 1$ we consider $K_x^k = ((f-t)^k,I_x)$.
Then $f(x)=c$ if and only if $K_x$ has non-zero dimension (as a vector space).
Moreover if $f(x) = c$ then, by the Nullstellensatz, $(f-t)^k$ is in $I_x$
for a sufficiently large $k$. For such a $k$, the dimension of $K_x$ is
the Milnor number at $x$ if $f(x)=c$, and it is $0$ otherwise.
Such a $k$ is less or equal to the Milnor number at $x$,
but $k$ can often be chosen much less. The minimal $k$ is the first integer
such that the vector space dimension of $K_x^k$ is equal to the one of  $K_x^{k+1}$.

\section{Milnor numbers and critical values at infinity}
\label{sec:inf}

We give the computation of the objects at infinity and its implementation 
in \Singular.
We will suppose that $f$ has isolated singularities at infinity, 
in fact computations are valid for a larger class of polynomials but
it cannot be computed if $f$ belongs to this class.
The algorithm is based on the article of D.~Siersma and M.~Tib\u ar,
\cite{ST}, that gives critical values at infinity and Milnor numbers
at infinity with the help of polar curves.

\subsection{Working space}

We will work in $\Pp^n \times \Cc$, with the homogeneous coordinates of $\Pp^n$:
$(x_1:\ldots:x_n:z)$ ; we still need $t$ which is a parameter or a variable depending
on the context.

We recall that
$$X = \left\lbrace ((x:z),t)\in \Pp^n\times \Cc \ | \ \bar
f(x,z)-tz^d=0 \right\rbrace.$$
The part at infinity of $X$ is $X_\infty = X \cap (\mathcal{H}_\infty \times \Cc)$:
$$X_\infty = \left\lbrace ((x:0),t)\in \Pp^n\times \Cc \ | \
f^d(x)=0 \right\rbrace,$$
Where $f = f^d + f^{d-1}+\cdots$ is the decomposition in
homogeneous polynomials.

In \Singular, we write:
\begin{verbatim}
    ring r = 0, (x(1..n),z,t), dp;
    poly f = ...;
    poly fH = homog(f,z)-t*z^deg(f);
    ideal X = fH;
    ideal Xinf = z, fH;
\end{verbatim}

\subsection{Polar curve}

Let $k$ be in $\{1,\ldots,n\}$.  The polar curve $\PP$
is the critical locus of the map $\phi : \Cc^n \longrightarrow \Cc^2$
defined for $x= (x_1,\ldots,x_n)$ by $\phi(x) = (f(x),x_k)$:
$$\PP = \left\lbrace
x \in \Cc^n \ | \ \frac{\partial f}{\partial x_i}(x) = 0, \forall i \not= k
\right\rbrace.$$

We have that $\PP$ is a curve or is void. 
We call $\PP_H$ the projective closure of $\PP$.
This curve intersects the hyperplane at infinity
$\mathcal{H}_\infty$ in finitely many points.
\begin{verbatim}
    ideal P = diff(f,x(1)),..., diff(f,x(k-1)), diff(f,x(k+1)),...;
    ideal PH = homog(P,z);
\end{verbatim}

The former objects can be viewed in $X$, we will also
denote by $\PP_H$, the set $(\PP_H \times \Cc) \cap X$.
In the chart $x_k=1$, we denote the curve $\PP_H$ by $\bar \CC$.
The ``real" polar curve $\CC$ in this chart is the closure of
$\bar \CC \setminus X_\infty$:
\begin{verbatim}
    ideal Cbar = x(k)-1, PH, X;
    ideal C = sat(Cbar,Xinf)[1];
\end{verbatim}

\subsection{Critical values at infinity}

We need the following result of \cite{ST}.
A value $c$ is a critical values at infinity
if and only there is coordinate $x_k$ and a point
$(x:0,t)$ in $X_\infty$ (with $x_k \not=0$) such that $(x:0,t) \in \CC$.
That is to say $\Binf$ is the projection of $\CC_\infty = X_\infty \cap \CC$ on
the space of parameters $t\in\Cc$.

Then the critical values are computed with:
\begin{verbatim}
    ideal Cinf = z, C;
    poly Binf = eliminate(Cinf,x(1)x(2)..x(n)z)[1];
\end{verbatim}

The set of critical values at infinity are the roots of the polynomial
\ttt{Binf}, which belongs to $\Cc[t]$.

\subsection{Milnor numbers at infinity}

Actually the results in \cite{ST} are more precise.
For a fixed $t$, let $X_t = \{ (x:z,t) \in X \}$, this is a projective model
for the fiber $f^{-1}(t)$.
\begin{theorem}[\cite{ST}]
The Milnor number at infinity at a point $(x:0,t) \in \CC_\infty$
is given by the intersection number (in $X$) of $\CC$ with $X_t$ at $(x:0,t)$.
\end{theorem}
So, for $c\in \Binf$, the Milnor number at infinity $\lambda_c$
(for the chart $x_k\not=0$),
is equal to the sum of all intersection numbers of $X_c$ and $\CC$
in $X_\infty$.

We compute an ideal $I$ which correspond to $X_c \cap \CC$,
then we only deals with points at infinity by intersecting it
this set with $z^q = 0$, for a sufficiently large $q$. 
\begin{verbatim}
    number c = ...;
    ideal Xc = t-c, X;
    ideal I = Xc, C;
    ideal K = z^q, I;    // q >> 1
    lambdac = vdim(std(K));
\end{verbatim}

Once we have computed $\lambda_c$ for all $c\in \Binf$, we have
$\lambda = \sum_{c\in \Binf}{\lambda_c}$.

\section{Families of polynomials}
\label{sec:fam}

Let $(f_s)_{s\in [0,1]}$ be a family of complex polynomials in $n$ variables.
We suppose that the coefficients are polynomial functions of $s$
and that for all $s\in [0,1]$, $f_s$ has affine isolated singularities 
and strong isolated singularities at infinity. The implementation is similar to the one
of paragraph \ref{sec:inf} and will be omitted.

\subsection{Change in affine space}

It is not possible to compute infinitely many $\mu(s)$, so we have to detect a change 
of $\mu(s)$. 
The Milnor numbers $\mu(s)$ changes if and only if some critical points
escape at infinity. Then we can detect critical parameters for $\mu$ as follows:
Let $J = \left\{ (x_1,\ldots, x_n,s) \in \Cc^n \times \Cc \mid 
\frac{\partial f_s}{\partial x_1}=\ldots, \frac{\partial f_s}{\partial x_n}=0\right\} $
be the set of critical points (that corresponds to the Jacobian ideal in $\Cc[x_1,\ldots,x_n,s]$).
Let $\bar J$ be the homogeneization of $J$ with the new variable $z$, while
$s$ is considered as a parameter. The part at infinity of $J$ corresponds
to the ideal $J_\infty = {\bar J} \cap (z=0)$, and the affine part of $J$ is
$J_\aff = \overline{\bar J \setminus J_\infty}$. Now the critical parameters for
$\mu$ is $\pr_s(J_\aff) \subset \Cc$, where $\pr_s$ is the projection to the $s$-coordinate.

It is possible to compute $\Baff(s)$ for all $s \in [0,1]$
by a direct extension of the work of paragraph \ref{sec:aff}.
Then we can compute the parameters where the cardinal of this set changes.

\subsection{Change at infinity}

Again it is not possible to compute infinitely many $\lambda(s)$.
We extend the definition of paragraph \ref{sec:inf} by adding a parameter $s$.
We set $d = \deg f_s$ and
$$X = \left\lbrace ((x:z),t,s)\in \Pp^n\times \Cc \times \Cc \ | \ \bar
f_s(x,z)-tz^d=0 \right\rbrace.$$
The part at infinity of $X$ is $X_\infty = X \cap (\mathcal{H}_\infty \times \Cc \times \Cc)$:
$$X_\infty = \left\lbrace ((x:0),t,s)\in \Pp^n\times \Cc \ | \
f_s^d(x)=0 \right\rbrace.$$
The polar ``curve'' is 
$$\PP = \left\lbrace
(x,s) \in \Cc^n \times \Cc \ | \ \frac{\partial f_s}{\partial x_i}(x) = 0, \forall i \not= k
\right\rbrace.$$
In the chart $x_k=1$ we denote the homogeneization of $\PP$ (with $s$ a parameter)
by $\bar{\mathcal{C}}$, and the ``real'' polar curve $\mathcal{C}$ in this chart is
the closure of $\bar{\mathcal{C}} \setminus X_\infty$.
The part at infinity of $\mathcal{C}$ is $\mathcal{C}_\infty = \mathcal{C} \cap X_\infty$.

Let $B_\infty(s)= \pr_t\{ (x:0,t,s) \in \mathcal{C}_\infty \}$.
For a generic $s'$, $\Binf(s') = B_\infty(s')$. 
 Then the critical parameters for $\Binf(s)$ is included in the set of parameters where
$\# B_\infty(s)$ fails to be equal to $\# \Binf(s')$ (in fact
$B_\infty(s)$ may be infinite).

We set $X_* = \{ (x:z,c,s) \in X \mid (x:0,c,s) \in \mathcal{C}_\infty \}$,
for non-critical parameters it corresponds to union of the irregular fibers at infinity.
Now a change of $\lambda$ corresponds a change in the value of the intersection multiplicity
of the polar curve $\mathcal{C}$ with $X_*$. 
The critical parameters for $\lambda$ are given as the projection to the $s$-coordinate
of $$ \overline{(\mathcal{C} \cap X_*) \setminus \mathcal{C}_\infty} \cap (z=0).$$

At last we compute parameters where the cardinal of $\B(s) = \Baff(s) \cup \Binf(s)$ changes.

\section{Examples}

\subsection{Brian\c{c}on polynomial}

The following example shows how to use the program once you have
started \Singular. We have to load the library \ttt{critic.lib}, then we set the ring,
with $n+1$ variables, the last variable will able to have the critical
values (as the zeroes of a polynomial) in return.
The following code gives critical values at infinity
of Brian\c{c}on polynomial.

\begin{verbatim}
    LIB "critic.lib";
    ring r = 0, (x,y,t), dp;
    poly s = xy+1;
    poly p = x*s+1;
    poly f = 3*y*p^3+3*p^2*s-5*p*s-s;
    crit(f);
\end{verbatim}
The result is:
\begin{verbatim}
    > Affine critical values are the roots of 1
    > Affine Milnor number : 0
    > Critical values at infinity are the roots of 3t2+16t
    > Milnor number at infinity : 4
    > Details of critical values at infinity :
    >   t       1
    >   3t+16   3
\end{verbatim}
This shows, that there is no affine critical value
(as the root of the polynomial $1$)
and that $\Binf = \{0,-\frac{16}{3}\}$ (as the root of the polynomial $t$ and $3t+16$)
are the critical values at infinity, with Milnor
number at infinity respectively equal to $1$ and $3$.

\subsection{More variables}

Let $f(a,b,c,d)= a+a^4b+b^2c^3+d^5$ be the example of
Choudary-Dimca, \cite{CD} and \cite{ACD}. This polynomial has isolated singularities
at infinity. The only singularity is a singularity at
infinity for the critical value $0$. Let's check it.

\begin{verbatim}
    ring r = 0, (a,b,c,d,t), dp;
    poly f = a+a^4*b+b^2*c^3+d^5;
    crit(f);
    > Affine critical values are the roots of 1
    > Affine Milnor number : 0
    > Critical values at infinity are the roots of t
    > Milnor number at infinity : 8
\end{verbatim}

\subsection{A family}

We give example of deformation, we first need to load the library
\ttt{defpoly.lib}, then we set a ring in $n+1$ variables, 
where the last variable is the parameter of the deformation.
For instance we consider the deformation $f_s(x,y)=y(1-sx)(y-(s-1)x)$.
\begin{verbatim}
    LIB "defpol.lib";
    ring r = 0, (x,y,s), dp;
    poly f = y*(1-sx)*(y-(s-1)*x);
    parCrit(f);
    > Critical parameters are included in the roots of s2-s
\end{verbatim}
Then the critical parameters are $s=0$ and $s=1$.

\subsection{A trivial family}

Another deformation is $f_s(x,y)= x(x^3y+sx^2+s^2x+1)$.
\begin{verbatim}
    LIB "defpol.lib";
    ring r = 0, (x,y,s), dp;
    poly f = x*(x^3*y+s*x^2+s^2*x+1);
    parCrit(f);
    > Critical parameters are included in the roots of 1
\end{verbatim}
Then $\mm(f_s)$ and the degree are constant; by Theorem \ref{th:mucst}
it implies that for all $s,s'\in \Cc$, $f_s$ and $f_{s'}$ are topologically equivalent.

\subsection{Combination}

We consider the family $f_s(x,y) = (x-s^2-1)(x^2y+1)$.
\begin{verbatim}
    LIB "defpol.lib";
    ring r = 0, (x,y,s), dp;
    poly f = (x-s^2-1)*(x^2*y+1);
    parCrit(f);
    > Critical parameters are included in the roots of s2+1
\end{verbatim}

For a generic value we have
\begin{verbatim}
    LIB "critic.lib";
    ring r = (0,s), (x,y,t), dp;
    poly f = (x-s^2-1)*(x^2*y+1);
    crit(f);
    > Affine critical values are the roots of t
    > Affine Milnor number : 1
    > Critical values at infinity are the roots of t+(s2+1)
    > Milnor number at infinity : 1
\end{verbatim}

And for a critical parameter ($s=i$ or $s=-i$):
\begin{verbatim}
    ring r = (0,s), (x,y,t), dp;
    minpoly = s^2+1;
    poly f = (x-s^2-1)*(x^2*y+1);
    crit(f);    
    > Affine critical values are the roots of 1
    > Affine Milnor number : 0
    > Critical values at infinity are the roots of t
    > Milnor number at infinity : 1
\end{verbatim}

\bigskip

{

}

\end{document}